\documentclass[leqno,draft,11pt,a4paper]{article}
\usepackage{amssymb,amsmath,theorem}
\usepackage[includeheadfoot,margin=1in]{geometry}

\allowdisplaybreaks[2]

\newcommand\eq{\leftrightarrow}
\newcommand\LOR{\bigvee}
\newcommand\ET{\bigwedge}
\newcommand\model{\vDash}

\newcommand\fii{\varphi}

\newcommand\p[1]{\langle#1\rangle}
\newcommand\lh[1]{\lvert#1\rvert}
\newcommand\bez{\smallsetminus}
\newcommand\sset{\subseteq}
\newcommand\nsset{\nsubseteq}
\newcommand\cupd{\mathbin{\dot\cup}}
\newcommand\nul{\varnothing}
\newcommand\res{\mathbin\restriction}

\DeclareMathOperator\dom{dom}
\DeclareMathOperator\Lh{lh}
\DeclareMathOperator\Th{Th}
\DeclareMathOperator\tcl{tc}
\DeclareMathOperator\stcl{\strc{tc}}
\DeclareMathOperator\rk{rk}

\newcommand\cxt[1]{\mathrm{#1}}
\newcommand\ptime{\cxt P}
\newcommand\nc{\cxt{NC}}
\newcommand\nci{\nc^1}
\newcommand\psp{\cxt{PSPACE}}
\newcommand\LT{\cxt{DLOGTIME}}
\newcommand\sksat[1]{S_{#1}\text-\mathrm{Sat}}
\newcommand\bsksat[1]{S_{#1}\text-\mathrm{BSat}}

\newcommand\vs{\!\mathit{VS}}
\newcommand\ax[1]{\mathrm{#1}}

\newcommand\N{\mathbb N}
\newcommand\strc[1]{\mathbf{#1}}
\newcommand\sH{\strc H}
\newcommand\sA{\strc A}
\newcommand\sB{\strc B}
\newcommand\sT{\strc T}
\newcommand\sM{\strc M}
\newcommand\sN{\strc N}
\newcommand\emb{\mathrel{\widetilde\sset}}
\newcommand\nemb{\mathrel{\widetilde\nsset}}

\newcommand\ob[1]{\overline{#1}}
\newcommand\txto{${}\to{}$}

\newenvironment{algo}[1][20em]{\catcode`\^^I=13 \obeylines\doalgo{#1}}{}
{\catcode`\^^I=13
\gdef\doalgo#1#2\end#{\hbox to\hsize{\hss \let^^I\qquad
  \def\>{{\everypar{}\indent}}%
  \fboxsep1em \linenum0 \parindent1.2em%
  \fbox{\hsize#1\vbox{%
  \everypar{\advance\linenum1 \llap{$\scriptstyle\the\linenum$\hskip.6em}}%
  #2}}\hss}\end}}
\newcount\linenum
\newcommand\key{\relax\ifmmode\expandafter\mathbf\else\expandafter\textbf\fi}

\newcommand\bme{\hskip.75em\relax}
\newcommand\noproof{\leavevmode\unskip\bme\vadjust{}\nobreak\hfill$\Box$\par}
\newenvironment{Pf}[1][]
  {\par\noindent\textit{Proof:}\bme\ignorespaces}
  {\noproof\pagebreak[2]\vskip\medskipamount\ignorespacesafterend}

\theoremstyle{plain}
\newtheorem{Thm}{Theorem}[section]
\newtheorem{Prop}[Thm]{Proposition}
\newtheorem{Cor}[Thm]{Corollary}
\newtheorem{Lem}[Thm]{Lemma}
\theorembodyfont\upshape
\newtheorem{Def}[Thm]{Definition}
\newtheorem{Rem}[Thm]{Remark}

\author{Emil Je\v r\'abek\thanks{Supported by grant 19-05497S of GA \v CR.
The Institute of Mathematics of the Czech Academy of Sciences is supported by RVO: 67985840.}\\[\medskipamount]
Institute of Mathematics, Czech Academy of Sciences\\
\small \v Zitn\'a 25,
115\:67 Praha 1,
Czech Republic,
email: \texttt{jerabek@math.cas.cz}
}

\title{The theory of hereditarily bounded sets}

\begin{document}
\maketitle

\begin{abstract}
We show that for any $k\in\omega$, the structure $\p{H_k,{\in}}$ of sets that are hereditarily of size at most~$k$ is
decidable. We provide a transparent complete axiomatization of its theory, a quantifier elimination result, and tight
bounds on its computational complexity. This stands in stark contrast to the structure
$V_\omega=\bigcup_kH_k$ of hereditarily finite sets, which is well known to be bi-interpretable with the
standard model of arithmetic $\p{\N,+,\cdot}$.
\end{abstract}

\section{Introduction}\label{sec:introduction}

The \emph{Vaught set theory $\vs$}, originally introduced by Vaught~\cite{vaught}, is a very rudimentary theory of
sets: it is axiomatized by the schema
\[\tag{$\ax V_n$}\forall x_0,\dots,x_{n-1}\:\exists y\:\forall t\:\Bigl(t\in y\eq\LOR_{i<n}t=x_i\Bigr)\]
for all $n\in\omega$, asserting that $\{x_i:i<n\}$ exists. It is one of the weakest known essentially undecidable
theories; while Robinson's theory~$R$, introduced in~\cite{tmr}, is even weaker (in terms of interpretability), $\vs$
is appealing in the simplicity of its axioms, especially in the context of set theories where setting up an
interpretation of an arithmetic theory such as~$R$ may be somewhat of a laborious task.

In contrast to full~$\vs$, the finite fragments $\vs_k$ (axiomatized by $(\ax V_0)$ and~$(\ax V_k)$, which imply
$(\ax V_m)$ for all $m\le k$) are \emph{not} essentially undecidable, but the reason for this is a bit indirect: for
each~$k$, $\vs_k$ is interpretable in any theory with a pairing function, and it is known that there exist decidable
consistent theories with pairing.

The first such theories were constructed by Malcev \cite{malc:loc-free,malc:ax-cl-loc-free}: he proved the decidability
of theories of \emph{locally free algebras}, which are essentially the first-order theories of term algebras in a given
signature. His results also apply to free algebras with function symbols constrained to be symmetric w.r.t.\ prescribed
groups of permutations of the arguments. As a special case, \emph{acyclic pairing functions} are locally free algebras
with a single binary function; e.g., the pairing function $2^x3^y$ on~$\N$ is acyclic, hence $\p{\N,2^x3^y}$ is
decidable. More generally, Tenney~\cite{tenn:phd} proved that pairing functions that are acyclic up to a finite (or
sufficiently well-behaved) set of exceptions have a decidable theory, including common pairing functions on~$\N$ such
as $2^x(2y+1)-1$, $\max\{x^2,y^2+x\}+y$, or Cantor's function $C(x,y)=\binom{x+y+1}2+x$. (The decidability
of $\p{\N,C}$ was reproved in~\cite{ceg-gri-rich:can-pair} using Malcev's results.) Decidable structures with pairing
may include more arithmetic functions: C\'egielski and Richard observed in~\cite{ceg-rich:lists} that pairing
functions such as $2^x+2^{x+y}$ are definable in $\p{\N,{+},2^x}$, which is decidable due to
Semenov~\cite{sem:plus-pow}, and in a tour de force~\cite{ceg-rich:can-pair-succ}, they proved the decidability of
$\p{\N,C,S}$ (while other related structures, including $\p{\N,C,{<}}$, $\p{\N,C,{+}}$, and $\p{\N,C,{\cdot}}$, are
undecidable).

For more background on theories with ``containers'' such as pairs, sets, and sequences, see Visser~\cite{vis:pair}.

While the results above confirm that finite fragments of the Vaught set theory are not essentially undecidable, the
decidable extensions of~$\vs_k$ we get from interpretation in theories of pairing are quite unnatural when we think of
them as set theories: for example, they will contradict extensionality, which is arguably the most characteristic
principle distinguishing sets from other kinds of objects. Thus, it might be worthwhile to see if we can find
decidable extensions of~$\vs_k$ that are easier to understand.

One of the simplest---and perhaps most natural---models of~$\vs_k$ is the structure%
\footnote{The standard notation in set theory is that, for a (usually regular infinite) cardinal~$\kappa$, $H_\kappa$
consists of sets hereditarily of cardinality ${<}\kappa$, thus our~$H_k$ would be denoted $H_{k+1}$. We decided to
violate this convention as it seems to be more confusing than helpful in the finite case.}
$\sH_k=\p{H_k,{\in}}$ of \emph{sets hereditarily of size at most~$k$}; that is, $H_k$ is the smallest family of sets
such that every subset of~$H_k$ of cardinality~${\le}k$ is a member of~$H_k$:
\[\forall x\:\bigl(x\sset H_k\land\lh x\le k\implies x\in H_k\bigr).\]
Equivalently, $H_k$ consists of (well-founded) sets $x$ such that $x$~itself, and all elements of its transitive
closure, have cardinality~${\le}k$.
The better known family of \emph{hereditarily finite} sets~$V_\omega$ includes each~$H_k$, and in fact,
$V_\omega=\bigcup_{k\in\omega}H_k$. We observe that $\sH_k$ is a \emph{minimal} model of~$\vs_k$, in that it
embeds (as a transitive submodel) into any other model of~$\vs_k$; thus, $\sH_k$ is canonically associated
with~$\vs_k$.

The main purpose of this paper is to show that $\Th(\sH_k)$ is decidable, providing an explicit natural example of a
decidable extension of~$\vs_k$. We present a transparent recursive axiomatization of $\Th(\sH_k)$, and a
characterization of elementary equivalence of tuples in models of $\Th(\sH_k)$ in terms of isomorphism of transitive
closures. Apart from the decidability of~$\sH_k$, this yields a quantifier elimination result (every formula is
equivalent to a Boolean combination of bounded existential formulas). We also establish that $\Th(\sH_k)$ is
stable, and it is not finitely axiomatizable. In Section~\ref{sec:comp-compl}, we investigate in more detail the
computational complexity of $\Th(\sH_k)$: we give an algorithm deciding $\Th(\sH_k)$ whose running time closely
matches a general lower bound on the complexity of theories with pairing by Ferrante and Rackoff~\cite{ferr-rack},
and its variant that has much lower complexity for sentences with a small number of quantifier alternations.

The properties of $\sH_k$ may be contrasted with the structure $\p{V_\omega,{\in}}$, which is bi-interpretable with
$\p{\N,{+},{\cdot}}$ by Ackermann~\cite{ack:hered-fin}, and as such it is heavily undecidable, and its quantifier
alternation hierarchy is proper.

Let us remark that while we formulate most results so that they apply to all $k\in\omega$, the cases $k=0,1$ are
somewhat degenerate: $\sH_0$ is a one-element structure, and $\sH_1$ is definitionally equivalent to
$\p{H_1,\nul,\{-\}}\simeq\p{\N,0,S}$. Moreover, the case $k=2$ can be reduced to Malcev's results: $\sH_2$ is
definitionally equivalent to the structure $\p{H_2,\nul,\{-,-\}}$, which is a free algebra with a constant and a
commutative binary operation. A similar reduction does not seem possible for $k\ge3$, as the set builder
operation $\{x_0,\dots,x_{k-1}\}$ has peculiar symmetries such as $\{x,x,y\}=\{x,y,y\}$ that cannot be expressed by mere
permutations of arguments.

\section{Completeness and decidability}\label{sec:compl-decid}

The main result of this section is the decidability of~$\sH_k$. Our strategy is to propose a recursively axiomatized
theory~$S_k$, true in~$\sH_k$, and prove its completeness: this implies that $S_k$ is decidable and $S_k=\Th(\sH_k)$.
Without further ado, here is the definition of~$S_k$.
\begin{Def}\label{def:sk}
Let $k\in\omega$. The theory $S_k$ in the language of set theory $\p{{\in}}$ is axiomatized by $(\ax V_0)$,
$(\ax V_k)$, the extensionality axiom
\[\tag{$\ax E$} \forall x,y\:\bigl(\forall t\:(t\in x\eq t\in y)\to x=y\bigr),\]
the boundedness axiom
\[\tag{$\ax B_k$} \forall x,u_0,\dots,u_k\:\Bigl(\ET_{i\le k}u_i\in x\to\LOR_{i<j\le k}u_i=u_j\Bigr)\]
postulating that all sets have at most $k$ elements, and the axioms
\[\tag{$\ax C_n$} \forall x_0,\dots,x_n\:\neg\Bigl(\ET_{i<n}x_i\in x_{i+1}\land x_n=x_0\Bigr)\]
for all $n\in\omega$, $n\ge1$, prohibiting finite $\in$-cycles.
\end{Def}

Clearly, $\sH_k\model S_k$. We aim to show that $S_k$ is complete; we will prove this by an Ehrenfeucht--Fra\"iss\'e
argument, which will more generally provide a characterization of elementary equivalence of finite tuples in models
of~$S_k$. Let us first agree on basic notation concerning models.
\begin{Def}\label{def:models}
As a general notational convention, we will denote first-order structures by bold-face letters (possibly decorated).
The domain of a structure will be denoted by the same letter, but in italics, and the basic relations and functions of
a structure carry the name of the structure as a superscript (this convention will also extend on a case-by-case basis
to various defined concepts). For example, a typical model of the language of set theory will be denoted~$\sA$, in
which case $\sA=\p{A,{\in}^\sA}$.

We denote finite tuples (sequences) by letters with bars such as $\ob a$; then $\Lh(\ob a)$ denotes the length
of~$\ob a$, and the individual elements of~$\ob a$ are $a_i$ with $0\le i<\Lh(\ob a)$.

Let $\sA$ and $\sB$ be structures for the same language, and $\ob a\in A$, $\ob b\in B$ finite
tuples of the same length $l=\Lh(\ob a)=\Lh(\ob b)$. We write $\sA,\ob a\equiv\sB,\ob b$ if $\ob a$ and $\ob b$
satisfy the same formulas, and $\sA,\ob a\equiv_n\sB,\ob b$ if they satisfy the same formulas of quantifier
rank at most~$n$. We recall that the quantifier rank of a formula~$\fii$ is defined inductively by
\begin{align*}
\rk(\fii)&=0,&&\text{$\fii$ quantifier-free,}\\
\rk\bigl(c(\fii_0,\dots,\fii_{k-1})\bigr)&=\max\bigl\{\rk(\fii_i):i<k\bigr\},&&c\in\{\land,\lor,\to,\neg\},\\
\rk(Q x\,\fii)&=\rk(\fii)+1,&&Q\in\{\exists,\forall\}.
\end{align*}

If $f\colon A\to B$ and $X\sset A$, then $f[X]$ denotes the image $\{f(x):x\in X\}$, and $f\res X\colon X\to B$ the
restriction of $f$ to~$X$. If $\ob a\in A$ with $l=\Lh(\ob a)$, then $f(\ob a)$ denotes the $l$-tuple $\ob b$ such that
$b_i=f(a_i)$ for each $i<l$.
\end{Def}

We also fix some notation and terminology specific to models of~$S_k$. In particular, we intend to characterize the
elementary equivalence relations $\equiv_n$ in terms of isomorphism of levels of transitive closures, hence we need to
define the latter.
\begin{Def}\label{def:tcequiv}
\emph{Bounded quantifiers} in the language of set theory are introduced as the abbreviations
\begin{align*}
\exists y\in x\:\fii&\equiv\exists y\:(y\in x\land\fii),\\
\forall y\in x\:\fii&\equiv\forall y\:(y\in x\to\fii),
\end{align*}
where $x$ and~$y$ are distinct variables. A formula is \emph{bounded} if it is built from atomic formulas using
Boolean connectives and bounded quantifiers.

If $\sA\model S_k$ and $u\in A$, then
\[u^\sA=\{v\in A:v\in^\sA u\}\]
denotes the \emph{extension} of $u$ in~$\sA$. Conversely, if $r\le k$ and $\{u_i:i<r\}\sset A$, then $\{u_i:i<r\}^\sA$
or $\{u_0,\dots,u_{r-1}\}^\sA$ denotes the $v\in A$ such that $v^\sA=\{u_i:i<r\}$, which exists by $(\ax V_0)$ or~$(\ax
V_k)$, and is unique by~$(\ax E)$. In particular, $\nul^\sA=\{\}^\sA$.

If $\ob a\in A$ and $l=\Lh(\ob a)$, we define levels of the \emph{transitive closure} of $\ob a$
(as subsets of $A$) by
\begin{align*}
\tcl^\sA_0(\ob a)&=\{a_i:i<l\},\\
\tcl^\sA_{n+1}(\ob a)&=\tcl^\sA_n(\ob a)\cup\bigcup_{u\in\tcl^\sA_n(\ob a)}u^\sA,\\
\tcl^\sA(\ob a)&=\bigcup_{n\in\omega}\tcl^\sA_n(\ob a).
\end{align*}
We denote by $\stcl^\sA_n(\ob a)$ the (possibly empty) structure $\p{\tcl^\sA_n(\ob a),{\in}^\sA,\ob a}$, and likewise,
$\stcl^\sA(\ob a)=\p{\tcl^\sA(\ob a),{\in}^\sA,\ob a}$.

Notice that $\tcl^\sA_n(\ob a)=\bigcup_{i<l}\tcl^\sA_n(a_i)$, and $\tcl^\sA_n(\ob a)$ is finite:
$\lh{\tcl^\sA_n(\ob a)}\le l\,k^{\le n}$, where
\[k^{\le n}=\sum_{i=0}^nk^i=\begin{cases}\dfrac{k^{n+1}-1}{k-1},&k\ne1,\\n+1,&k=1.\end{cases}\]
Also, for any fixed $n$ and~$l$, there is a formula $\fii(\ob x,y)$ with $\Lh(\ob x)=l$ that defines the relation
$y\in\tcl^\sA_n(\ob x)$ in every model $\sA\model S_k$. Finally, we define
\begin{align*}
\sA,\ob a\sim_n\sB,\ob b&\iff\stcl^\sA_n(\ob a)\simeq\stcl^B_n(\ob b),\\
\sA,\ob a\sim_{\phantom n}\sB,\ob b&\iff\stcl^A(\ob a)\simeq\stcl^B(\ob b).
\end{align*}
\end{Def}

We first observe basic properties of morphisms on transitive closures.
\begin{Lem}\label{lem:tc-restr}
Let $\sA,\sB\model S_k$, $\ob a\in A$, $\ob b\in B$, $\Lh(\ob a)=\Lh(\ob b)$, and $n>0$.
\begin{enumerate}
\item\label{item:3}
If $f\colon\tcl_n^\sA(\ob a)\to B$ is a mapping such that $f(\ob a)=\ob b$ and
\begin{equation}\label{eq:3}
\forall u\in\tcl_{n-1}^\sA(\ob a)\:f(u)=\{f(t):t\in^\sA u\}^\sB,
\end{equation}
then $f[\tcl_m^\sA(\ob a)]=\tcl_m^\sB(\ob b)$ for all $m\le n$.
\item\label{item:4}
Any $f\colon\stcl_n^\sA(\ob a)\simeq\stcl_n^\sB(\ob b)$ satisfies~\eqref{eq:3}, thus
$f\res\tcl_m^\sA(\ob a)\colon\stcl_m^\sA(\ob a)\simeq\stcl_m^\sB(\ob b)$ for all $m\le n$.
\end{enumerate}
\end{Lem}
\begin{Pf}

\ref{item:3}: By induction on~$m$. The case $m=0$ holds. Assume $f[\tcl_m^\sA(\ob a)]=\tcl_m^\sB(\ob b)$ and $m<n$. If
$t\in\tcl_{m+1}^\sA(\ob a)$, then $t\in^\sA u$ for some $u\in\tcl_m^\sA(\ob a)$, thus $f(t)\in^\sB
f(u)\in\tcl_m^\sB(\ob b)$ by \eqref{eq:3} and the induction hypothesis, which means $f(t)\in\tcl_{m+1}^\sB(\ob b)$.
Conversely, if $s\in\tcl_{m+1}^\sB(\ob b)$, we have $s\in^\sB v$ for some $v\in\tcl_m^\sB(\ob b)$. By the induction
hypothesis, there is $u\in\tcl_m^\sA(\ob a)$ such that $f(u)=v$, thus $s=f(t)$ for some $t\in^\sA u$ by~\eqref{eq:3},
whence $t\in\tcl_{m+1}^\sA(\ob a)$.

\ref{item:4}: Let $u\in\tcl_{n-1}^\sA(\ob a)$. We can prove $f(u)\in\tcl_{n-1}^\sB(\ob b)$ as in~\ref{item:3}. On the one
hand, if $t\in^\sA u$, then $t\in\dom(f)$, hence $f(t)\in^\sB f(u)$. Now, on the other hand, if $s\in^\sB f(u)$, then
$s\in\tcl_n^\sB(\ob b)$, which means that $s=f(t)$ for some~$t\in\tcl_n^\sA(\ob a)$. Then $f(t)\in^\sB f(u)$ implies
$t\in^\sA u$.
\end{Pf}

By definition, ${\equiv}=\bigcap_n{\equiv}_n$. It may not be a priori obvious that the same holds for the $\sim$
relation (which we aim to eventually prove to coincide with~$\equiv$): e.g., the corresponding property fails for
general pointed directed acyclic graphs. However, here it is true because axiom $(\ax B_k)$ ensures that the graphs are
image-finite:
\begin{Lem}\label{lem:sim-simn}
Let $\sA,\sB\model S_k$, $\ob a\in A$, $\ob b\in B$, and $\Lh(\ob a)=\Lh(\ob b)$. Then $\sA,\ob a\sim\sB,\ob b$ if and
only if $\forall n\in\omega\,\sA,\ob a\sim_n\sB,\ob b$.
\end{Lem}
\begin{Pf}
The left-to-right implication is clear. For the converse, Lemma~\ref{lem:tc-restr} shows that the set~$T$ of all
isomorphisms $f\colon\stcl^\sA_n(\ob a)\simeq\stcl^\sB_n(\ob b)$, $n\in\omega$, forms a tree when ordered by inclusion,
and the finiteness of $\tcl_n$ implies that $T$ is finitely branching. As such, $T$ has an infinite branch by K\H
onig's lemma; the union of the branch is then an isomorphism of $\stcl^\sA(\ob a)$ to $\stcl^\sB(\ob b)$.
\end{Pf}

It is relatively straightforward to prove that $\sA,\ob a\equiv\sB,\ob b$ implies $\sA,\ob a\sim\sB,\ob b$: in view of
the previous lemma, we only need to establish that the isomorphism types of the finite structures $\stcl^\sA_n(\ob a)$
are definable. We do this below, including explicit bounds on the complexity of the defining formulas.
\begin{Lem}\label{lem:tc-fla}
Let $\sA\model S_k$, $\ob a\in A$, $l=\Lh(\ob a)$, and $n\in\omega$. Then there is a formula $\fii_{\ob a,n}(\ob x)$ such
that for any $\sB\model S_k$ and any $l$-tuple $\ob b\in B$, we have
\[\sB\model\fii_{\ob a,n}(\ob b)\iff\sA,\ob a\sim_n\sB,\ob b.\]
Moreover, we may take $\fii_{\ob a,n}$ in the form $\psi(\ob x)\land\neg\LOR_{i<m}\psi_i(\ob x)$, where $\psi$
and~$\psi_i$ are bounded existential formulas using at most $l(k^{\le n}-1)$ quantifiers each.
\end{Lem}
\begin{Pf}
Let $\{a_i:l\le i<r\}$ be an enumeration of $\tcl^\sA_n(\ob a)\bez\{a_i:i<l\}$, where $r\le l\,k^{\le n}$, and for every
$i\ge l$, there is $p(i)<i$ such that $a_i\in^\sA a_{p(i)}$ (this can be arranged by enumerating elements of
$\tcl^\sA_{n'}(\ob a)$ before elements of $\tcl^\sA_{n'+1}(\ob a)\bez\tcl^\sA_{n'}(\ob a)$, for each $n'<n$). Let
$\theta$ be (the conjunction of) the diagram of $\{a_i:i<r\}$ with the structure induced from~$\sA$, and put
\[\psi(x_0,\dots,x_{l-1})=
 \exists x_l\in x_{p(l)}\:\exists x_{l+1}\in x_{p(l+1)}\:\dots\:\exists x_{r-1}\in x_{p(r-1)}\:\theta(x_0,\dots,x_{r-1}).\]
Then for any $\sB\model S_k$ and $\ob b\in B$,
\[\sB\model\psi(\ob b)\iff\stcl^\sA_n(\ob a)\emb\stcl^\sB_n(\ob b),\]
where $\sM\emb\sN$ denotes that there exists an embedding $f\colon\sM\to\sN$. Let $\{\sM_i:i<m\}$
be an enumeration (up to isomorphism) of all structures of the form $\stcl^\strc C_n(\ob c)$ that do not embed into
$\stcl^\sA_n(\ob a)$, and as above, let $\psi_i$ be a bounded existential formula in at most $l\,k^{\le n}$ variables
such that
\[\sB\model\psi_i(\ob b)\iff\sM_i\emb\stcl^\sB_n(\ob b).\]
Then $\fii_{\ob a,n}=\psi\land\neg\LOR_{i<m}\psi_i$ satisfies
\begin{align*}
\sB\model\fii_{\ob a,n}(\ob b)
&\iff\stcl^\sA_n(\ob a)\emb\stcl^\sB_n(\ob b)\land\forall i<m\:\sM_i\nemb\stcl^\sB_n(\ob b)\\
&\iff\stcl^\sA_n(\ob a)\emb\stcl^\sB_n(\ob b)\land\stcl^\sB_n(\ob b)\emb\stcl^\sA_n(\ob a)\\
&\iff\stcl^\sA_n(\ob a)\simeq\stcl^\sB_n(\ob b),
\end{align*}
using the fact that if $\sM$ and~$\sN$ are finite structures such that $\sM\emb\sN$ and $\strc
N\emb\sM$, then $\sM\simeq\sN$.
\end{Pf}
\begin{Cor}\label{cor:eleq-sim}
Let $\sA,\sB\model S_k$, $\ob a\in A$, $\ob b\in B$, and $l=\Lh(\ob a)=\Lh(\ob b)$. Then $\sA,\ob a\equiv\sB,\ob b$
implies $\sA,\ob a\sim\sB,\ob b$. More precisely, $\sA,\ob a\equiv_{l(k^{\le n}-1)}\sB,\ob b$ implies $\sA,\ob
a\sim_n\sB,\ob b$.
\noproof\end{Cor}

It is more difficult to show the converse implication $\sA,\ob a\sim\sB,\ob b\implies\sA,\ob a\equiv\sB,\ob b$. We will
do it by an Ehrenfeucht--Fra\"\i ss\'e argument: that is, we will prove that if $\sA,\ob a\sim_m\sB,\ob b$ for $m$
sufficiently larger than~$n$, then any extension of $\ob a$ to $\ob a,c$ can be matched by an extension of $\ob b$ to
$\ob b,d$ so that $\sA,\ob a,c\sim_n\sB,\ob b,d$. This is the content of the crucial Lemma~\ref{lem:simn-ext} below.
However, we start with a little technical result that will be needed in its proof.
\begin{Lem}\label{lem:gen-pts}
Let $\sA\model S_k$, $\ob a\in A$, and $n,r\in\omega$, where $k\ge1$. There exists $\{v_i:i<r\}\sset A$ such that
\begin{itemize}
\item $v_i\notin\tcl^\sA(\ob a)$,
\item $i\ne j\implies v_i\notin\tcl^\sA_n(v_j)$,
\end{itemize}
for all $i,j<r$.
\end{Lem}
\begin{Pf}
We may assume that $l=\Lh(\ob a)>0$. The acyclicity of $\in^\sA$ implies that the relation $x\in\tcl^\sA(y)$ (which is
the reflexive transitive closure of $\in^\sA$) is a partial order, hence its restriction to any nonempty finite set has
a maximal element. That is, we can find $a\in\{a_i:i<l\}$ such that $a\notin\tcl^\sA(a_i)$ for any $a_i\ne a$. Then
$\{a\}^\sA\notin\tcl^\sA(\ob a)$, which implies that $v_i=\{a\}^{1+(n+1)i}$ have the required properties, where
$\{a\}^0=a$, $\{a\}^{t+1}=\{\{a\}^t\}^\sA$. (If $k\ge2$, we may even ensure the stronger condition
$v_i\notin\tcl^\sA(v_j)$ for $j\ne i$, by putting $v_i=\{\{a\}^{i+1},\{a\}^{i+2}\}^\sA$.)
\end{Pf}
\begin{Lem}\label{lem:simn-ext}
Let $\sA,\sB\model S_k$, $\ob a\in A$, $\ob b\in B$, $l=\Lh(\ob a)=\Lh(\ob b)$, and $n>0$. If
$\sA,\ob a\sim_{k^{\le n}+n}\sB,\ob b$, then for every $c\in A$, there exists $d\in B$ such that $\sA,\ob
a,c\sim_{n-1}\sB,\ob b,d$.
\end{Lem}
\begin{Pf}
If $k=0$, the conclusion of the lemma holds trivially as $\lh A=\lh B=1$, hence we may assume $k\ge1$. Put
$N=k^{\le n}+n$, and fix $f\colon\stcl^\sA_N(\ob a)\simeq\stcl^\sB_N(\ob b)$. Let $C$ be the smallest subset
of~$\tcl^\sA_n(c)\bez\tcl^\sA_n(\ob a)$ satisfying the inductive condition
\[u^\sA\sset\tcl^\sA_n(\ob a)\cup C\implies u\in C\]
for $u\in\tcl^\sA_n(c)\bez\tcl^\sA_n(\ob a)$. We can extend $f\res\tcl^\sA_n(\ob a)$ uniquely to a mapping
$g\colon\tcl^\sA_n(\ob a)\cup C\to B$ such that
\[g(u)=\{g(t):t\in^\sA u\}^\sB\]
for all $u\in C$. Let $\{u_i:i<r\}$ be an injective enumeration of
\[\{u\in\tcl^\sA_n(c)\bez\tcl^\sA_n(\ob a):u^\sA\nsset\tcl^\sA_n(\ob a,c)\}.\]
Using Lemma~\ref{lem:gen-pts}, we can find $\{v_i:i<r\}\sset B$ such that
\begin{enumerate}
\item\label{item:5} $v_i\notin\tcl^\sB_N(\ob b)\cup g[C]$,
\item\label{item:6} $i\ne j\implies v_i\notin\tcl^\sB_N(v_j)$, 
\end{enumerate}
for all $i,j<r$. Since $\in^\sA$ is acyclic, and therefore well-founded on the finite set $\tcl^\sA_n(c)$, we can
construct using well-founded recursion a unique mapping $g\colon\tcl^\sA_n(\ob a,c)\to B$ such that
\[g(u)=\begin{cases}
f(u),&u\in\tcl^\sA_n(\ob a),\\
v_i,&u=u_i,\\
\{g(t):t\in^\sA u\}^\sB,&u\in\tcl^\sA_n(c)\bez\tcl^\sA_n(\ob a),u^\sA\sset\tcl^\sA_n(\ob a,c).
\end{cases}\]
(This agrees with the original definition of~$g$ on $\tcl^\sA_n(\ob a)\cup C$, hence keeping the same name will not
lead to confusion. The reason for this slightly awkward two-stage construction of~$g$ is that we could not define the
whole $g$ right away as it depends on the choice of $\{v_i:i<r\}$, which in turn depends on $g\res C$.) Using
Lemma~\ref{lem:tc-restr} and the definition of~$g$, the condition
\begin{equation}\label{eq:4}
g(u)=\{g(t):t\in^\sA u\}^\sB
\end{equation}
holds for all $u\in\tcl^\sA_n(\ob a,c)$ such that $u^\sA\sset\tcl^\sA_n(\ob a,c)$. In particular, it holds for all
$u\in\tcl^\sA_{n-1}(\ob a,c)$, hence Lemma~\ref{lem:tc-restr} implies
$g[\tcl^\sA_{n-1}(\ob a,c)]=\tcl^\sB_{n-1}(\ob b,d)$, where $d=g(c)$.

We claim that $g$ is injective. Assuming for the moment that this is true, let us show that $g\res\tcl^\sA_{n-1}(\ob
a,c)\colon\stcl^\sA_{n-1}(\ob a,c)\simeq\stcl^\sB_{n-1}(\ob b,d)$. If $t,u\in\tcl^\sA_{n-1}(\ob a,c)$, then $u$
satisfies~\eqref{eq:4}. Thus, on the one hand, $t\in^\sA u$ implies $g(t)\in^\sB g(u)$; on the other hand, if
$g(t)\in^\sB g(u)$, then $g(t)=g(t')$ for some $t'\in^\sA u$, and we have $t=t'$ by injectivity, hence $t\in^\sA u$.

It remains to prove the injectivity of~$g$. Assume for contradiction that there are $x,y\in\tcl^\sA_n(\ob a,c)$ such
that $x\ne y$, but $g(x)=g(y)$. Since $\in^\sA$ is well-founded on $\tcl^\sA_n(\ob a,c)$, we may take $x$ to be
$\in^\sA$-minimal for which such a~$y$ exists.

If $x^\sA,y^\sA\sset\tcl^\sA_n(\ob a,c)$ so that \eqref{eq:4} holds for both $x$ and~$y$, there is $x'\in^\sA x$
such that $x'\notin^\sA y$, or $y'\in^\sA y$ such that $y'\notin^\sA x$. In the former case, $g(x')\in^\sB g(x)=g(y)$,
hence $g(x')=g(y')$ for some $y'\in^\sA y$, and necessarily $x'\ne y'$; this contradicts the minimality of~$x$. The
other case is symmetric.

Thus, $x^\sA\nsset\tcl^\sA_n(\ob a,c)$ or $y^\sA\nsset\tcl^\sA_n(\ob a,c)$. By swapping $x$ and~$y$ if necessary
(dropping the minimality assumption, which is no longer needed), we may assume the latter. We distinguish two cases.

\textbf{Case 1: $y=u_i$ for some $i<r$}.
Thus, $g(x)=v_i$ and $x\ne u_i$. We cannot have $x\in\tcl^\sA_n(\ob a)\cup C$
because of~\ref{item:5}, hence $x\in\tcl^\sA_n(c)\bez(\tcl^\sA_n(\ob a)\cup C)$. Put $x_0=x$. Either $x_0=u_j$ for
some~$j$, or $x_0^\sA\sset\tcl^\sA_n(\ob a,c)$, while $x_0^\sA\nsset\tcl^\sA_n(\ob a)\cup C$; thus, there is
$x_1\in^\sA x_0$ such that $x_1\in\tcl^\sA_n(c)\bez(\tcl^\sA_n(\ob a)\cup C)$, and we can continue in the
same way. By acyclicity of~$\in^\sA$, the process has to stop after less than $\lh{\tcl^\sA_n(c)}\le k^{\le n}$ steps;
that is, we can construct a sequence $x_0,\dots,x_s\in\tcl^\sA_n(c)\bez(\tcl^\sA_n(\ob a)\cup C)$ such that $s<k^{\le n}$,
$x_s\in^\sA\dots\in^\sA x_1\in^\sA x_0$, $x_i^\sA\sset\tcl^\sA_n(\ob a,c)$ for each $i<s$, and
$x_s^\sA\nsset\tcl^\sA_n(\ob a,c)$, which means $x_s=u_j$ for some~$j<r$. But then
\[v_j=g(x_s)\in^\sB\dots\in^\sB g(x_1)\in^\sB g(x_0)=v_i\]
by~\eqref{eq:4}, i.e., $v_j\in\tcl^\sB_s(v_i)$. By condition~\ref{item:6}, this is only possible if $j=i$, and then
$s=0$ by acyclicity of~$\in^\sB$. Thus, $x=u_i$ after all, a contradiction.

\textbf{Case 2: $y\in\tcl^\sA_n(\ob a)$}.
We cannot have $x\in\tcl^\sA_n(\ob a)$ as $f$ is injective. If $x\notin\tcl^\sA_n(\ob a)\cup C$, then the argument in
Case~1 shows that $v_j\in\tcl^\sB_{k^{\le n}}(g(x))$ for some $j<r$, while $g(x)=f(y)\in\tcl^\sB_n(\ob b)$, thus
$v_j\in\tcl^\sB_N(\ob b)$, contradicting~\ref{item:5}. The only remaining possibility is $x\in C$. Put $x_0=x$ and
$y_0=y$. We have $x_0^\sA\sset\tcl^\sA_n(\ob a)\cup C$, thus $x_0$ satisfies~\eqref{eq:4}, while
\[f(y_0)=\{f(t):t\in^\sA y_0\}^\sB\]
by Lemma~\ref{lem:tc-restr}. Thus, the same argument as above shows that there are $x_1\in^\sA x_0$ (whence
$x_1\in\tcl^\sA_n(\ob a)\cup C$) and $y_1\in^\sA y_0$ (whence $y_1\in\tcl^\sA_{n+1}(\ob a)$) such that $g(x_1)=f(y_1)$
and $x_1\ne y_1$. If $x_1\in C$, we may continue in the same way, but the acyclicity of $\in^\sA$ again implies that
the process has to stop: that is, we construct sequences $x_0,\dots,x_s$ and $y_0,\dots,y_s$ such that $s\le\lh C\le
k^{\le n}$, $x_s\in^\sA\dots\in^\sA x_1\in^\sA x_0$, $x_i\in C$ for each $i<s$, $x_s\in\tcl^\sA_n(\ob a)$,
$y_s\in^\sA\dots\in^\sA y_1\in^\sA y_0$ (thus $y_i\in\tcl^\sA_{n+i}(\ob a)\sset\tcl^\sA_N(\ob a)$), $x_i\ne y_i$ for
each $i\le s$, and $g(x_i)=f(y_i)$. But then $f(x_s)=f(y_s)$ contradicts the injectivity of~$f$. This completes the
proof.
\end{Pf}

We can now put everything together to obtain the desired characterization of elementary equivalence.
\begin{Thm}\label{thm:sim-eleq}
Let $\sA,\sB\model S_k$, $\ob a\in A$, $\ob b\in B$, and $l=\Lh(\ob a)=\Lh(\ob b)$. Then
\[\sA,\ob a\equiv\sB,\ob b\iff\sA,\ob a\sim\sB,\ob b.\]
More precisely, for all $n\in\omega$,
\begin{align}
\label{eq:1}\sA,\ob a\equiv_{l(k^{\le n}-1)}\sB,\ob b&\implies\sA,\ob a\sim_n\sB,\ob b,\\
\label{eq:2}\sA,\ob a\sim_{t_k(n)}\sB,\ob b&\implies\sA,\ob a\equiv_n\sB,\ob b,
\end{align}
where $t_k(0)=0$, $t_k(n+1)=k^{\le t_k(n)+1}+t_k(n)+1$.
\end{Thm}
\begin{Pf}
Corollary~\ref{cor:eleq-sim} gives~\eqref{eq:1}, hence it suffices to establish~\eqref{eq:2}. Clearly
\[\sA,\ob a\sim_{t_k(n)}\sB,\ob b\implies\sA,\ob a\equiv_0\sB,\ob b,\]
and by Lemma~\ref{lem:simn-ext} and the definition of $t_k(n+1)$, 
\begin{align*}
\sA,\ob a\sim_{t_k(n+1)}\sB,\ob b\implies{}&\forall c\in A\:\exists d\in B\:(\sA,\ob a,c\sim_{t_k(n)}\sB,\ob b,d)\\
{}\land{}&\forall d\in B\:\exists c\in A\:(\sA,\ob a,c\sim_{t_k(n)}\sB,\ob b,d).
\end{align*}
Thus, if $\sA,\ob a\sim_{t_k(n)}\sB,\ob b$, then Duplicator has a winning strategy in the $n$-round
Ehrenfeucht--Fra\"iss\'e game for $\p{\sA,\ob a}$ and $\p{\sB,\ob b}$, which implies $\sA,\ob a\equiv_n\sB,\ob b$.
\end{Pf}
\begin{Thm}\label{thm:dec-compl}
The theory $S_k$ is complete for each $k\in\omega$. Consequently, $S_k=\Th(\sH_k)$, and $S_k$ is decidable.
\end{Thm}
\begin{Pf}
Applying Theorem~\ref{thm:sim-eleq} with $l=0$, we see that any two models of~$S_k$ are elementarily equivalent, thus $S_k$
is complete. Being a complete recursively axiomatized theory, it is decidable.
\end{Pf}

In order to clarify the numerical content of Theorem~\ref{thm:sim-eleq}, let us give bounds on~$t_k$ using better known
functions.
\begin{Def}\label{def:supexp}
The iterated exponential function $2^x_n$ is defined by $2^x_0=x$ and $2^x_{n+1}=2^{2^x_n}$. Unless stated otherwise,
$\log n$ denotes logarithm to base~$2$.
\end{Def}
\begin{Prop}\label{prop:tk-growth}
We have $t_1(n)=3(2^n-1)$. For $k\ge2$ and $n\ge1$,
\begin{equation}\label{eq:5}
t_k(n)\le2^{c_k}_{n-1},
\end{equation}
where $c_k=(k+3)\log k+\log\log k+2$.
\end{Prop}
\begin{Pf}
The expression $t_1(n)=3(2^n-1)$ follows by induction on~$n$ from the defining recurrence, which simplifies to
$t_1(0)=0$, $t_1(n+1)=2t_1(n)+3$. For $k\ge2$, we put $f(x)=k^{\le x+1}+x+1$ and $h(x)=(x+1)\log k+\log\log k+2$. We
want to show
\begin{equation}\label{eq:6}
f^{(n)}(x)\le2^{h(x)}_n
\end{equation}
for all $x\ge0$ and $n\in\omega$, which gives~\eqref{eq:5} using $t_k(n)=f^{(n)}(0)=f^{(n-1)}(k+2)$.

Now, using the monotonicity of $2_n^x$ in~$x$, \eqref{eq:6} follows by induction on~$n$ from the inequalities
$x\le h(x)$ (which is obvious) and $h(f(x))\le2^{h(x)}$, hence it suffices to prove the latter; unwinding the
definitions, we need to show that
\begin{equation}\label{eq:7}
(k^{\le x+1}+x+2)\log k+\log\log k+2\le4k^{x+1}\log k=2^{(x+1)\log k+\log\log k+2}.
\end{equation}
This follows from the inequalities
\begin{gather*}
k^{\le x+1}\le\frac k{k-1}k^{x+1}\le2k^{x+1},\\
x+2\le2^{x+1}\le k^{x+1},\\
\log\log k+2\le2^{\log\log k+1}\le k\log k,
\end{gather*}
which are easy to verify, using the fact that $2^x\ge x+1$ for all $x\ge1$.
\end{Pf}
\begin{Rem}\label{rem:bounds-k=1}
For $k=1$, the bound $t_1(n)=3(2^n-1)$ from Theorem~\ref{thm:sim-eleq} can be improved to $2^n-1$, because in this case
Lemma~\ref{lem:simn-ext} holds with the conclusion strengthened to $\sA,\ob a,c\sim_n\sB,\ob b,d$. Moreover, one can also
prove a matching improvement to Corollary~\ref{cor:eleq-sim} to obtain the exact characterization
\[\sA,\ob a\equiv_n\sB,\ob b\iff\sA,\ob a\sim_{2^n-1}\sB,\ob b,\]
using the fact that there are definitions of quantifier rank~$n$ of $y=\{x\}^t$ for each $t\le2^n$ and of $y=\{\nul\}^t$
for each $t\le2^n-2$. We leave the details to an interested reader.
\end{Rem}

Apart from the completeness and decidability of~$S_k$, Theorem~\ref{thm:sim-eleq} implies a quantifier elimination result
for~$S_k$:
\begin{Thm}\label{thm:qe}
Let $k\in\omega$. Then every formula is equivalent to a Boolean combination of bounded existential formulas over~$S_k$.
\end{Thm}
\begin{Pf}
Let $\fii(\ob x)$ be a formula. By Theorem~\ref{thm:sim-eleq}, there exists $n$ such that
\[\sA,\ob a\sim_n\sB,\ob b\implies\bigl(\sA\model\fii(\ob a)\iff\sB\model\fii(\ob b)\bigr)\]
for any $\sA,\sB\model S_k$ and $\ob a\in A$, $\ob b\in B$. There are only finitely many isomorphism types of
structures of the form $\stcl_n^\sA(\ob a)$, thus there is a finite list $\{\p{\sA^i,\ob a^i}:i<m\}$ such that
\[\sA\model\fii(\ob a)\iff\exists i<m\:(\sA,\ob a\sim_n\sA^i,\ob a^i).\]
Then
\[S_k\vdash\fii(\ob x)\eq\LOR_{i<m}\fii_{\ob a^i,n}(\ob x),\]
where $\fii_{\ob a^i,n}(\ob x)$ is as in Lemma~\ref{lem:tc-fla}, which makes the right-hand side a Boolean combination of
bounded existential formulas.
\end{Pf}

\begin{Rem}\label{rem:ktuple-symbol}
If we expand the language with the predicates $y=\nul$ and $y=\{x_0,\dots,x_{k-1}\}$ (which have bounded universal
definitions in the original language), every formula is equivalent both to a bounded existential formula and to a
bounded universal formula. To see this, note that an embedding $f\colon\stcl_n^\sA(\ob a)\to\stcl_n^\sB(\ob b)$ in the
expanded language has to be an isomorphism, as $f[\tcl_n^\sA(\ob a)]=\tcl_n^\sB(\ob b)$ by Lemma~\ref{lem:tc-restr}. It
follows that if we take $\theta$ in the proof of Lemma~\ref{lem:tc-fla} to be the diagram in the expanded language, then
it suffices to put $\fii_{\ob a,n}=\psi$.

If $k=1$, the $y=\{x\}$ predicate is redundant, as it is equivalent to $x\in y$. Moreover, $S_1$ has full quantifier
elimination in a language with function symbols $\nul$ and~$\{x\}$, as $\p{H_1,\nul,\{x\}}\simeq\p{\N,0,S}$.
\end{Rem}

As we learned from Albert Visser, it is an interesting problem whether there exists a \emph{finitely axiomatized}
consistent decidable theory with a pairing function. We observe that our theories do not cut the mustard, though we
postpone the (albeit simple) proof to the next section, where the relevant construction will be used in a more
substantial way:
\begin{Prop}\label{prop:fin-ax}
$S_k$ is not finitely axiomatizable for any $k>0$.
\end{Prop}
\begin{Pf}
See Corollary~\ref{cor:fin-ax}.
\end{Pf}

\begin{Rem}\label{rem:found}
The axioms $(\ax C_n)$ of $S_k$ express the acyclicity of~$\in$. More generally, since $\in$ is well founded, $\sH_k$
satisfies the \emph{$\in$-induction} schema
\[\forall x\:\bigl(\forall y\in x\:\fii(y)\to\fii(x)\bigr)\to\forall x\:\fii(x)\]
(where $\fii$ is any formula, possibly with parameters), of which each $(\ax C_n)$ axiom is a special case. By
Theorem~\ref{thm:dec-compl}, the full $\in$-induction schema is equivalent to its instances $\{(\ax C_n):n\ge1\}$ over
the remaining axioms of~$S_k$; there does not seem to be an easy direct proof of this fact.

The axiom of foundation (regularity) as commonly formulated in ZF,
\[x\ne\nul\to\exists y\in x\:\forall z\:\neg(z\in x\land z\in y),\]
is strictly weaker: using the fact that $\lh x\le k$ by $(\ax B_k)$, it is easily seen to be equivalent to $\{(\ax
C_n):1\le n\le k\}$.
\end{Rem}

We end this section with a basic model-theoretic classification of the $S_k$~theories.
\begin{Def}\label{def:stab}
Let $\kappa$ be an infinite cardinal. A theory $T$ is \emph{$\kappa$-stable} if for every $\sM\model T$ and $A\sset M$
of size $\lh A\le\kappa$, there are at most~$\kappa$ (complete) $1$-types of~$\sM$ over~$A$. We say that $T$ is
\emph{stable} if it is $\kappa$-stable for some~$\kappa\ge\|T\|$, and it is \emph{superstable} if there is $\kappa_0$
such that $T$ is $\kappa$-stable for all $\kappa\ge\kappa_0$.
\end{Def}

As is well known, the theory $S_1$---definitionally equivalent to $\Th(\N,0,S)$---is uncountably categorical, and
therefore $\kappa$-stable for all $\kappa\ge\omega$. In contrast to that, it is easy to see that no consistent theory
with pairing (even non-functional) can be superstable, as there are always at least $\lh A^\omega$ different types
over~$A$; thus, the result below is the best possible for $k\ge2$.
\pagebreak[2]
\begin{Prop}\label{prop:stab}
For each $k\ge2$, the theory $S_k$ is stable.
\end{Prop}
\begin{Pf}
Let $\sM\model S_k$, and $A\sset M$ be such that $\lh A\le\kappa$. By replacing $A$ with $\tcl^\sM(A)$ (which has the
same cardinality) if necessary, we may assume $\tcl^\sM(A)=A$. By Theorem~\ref{thm:sim-eleq}, $1$-types over~$A$ correspond
to isomorphism types of $\stcl^\sN(A,c)$ for $\sN\succeq\sM$, $c\in N$; since the structure on~$A$ is fixed, these are
determined by isomorphism types of $\stcl^\sN(c)$ expanded with constants $a$ for all $a\in A\cap\tcl^\sN(c)$. In other
words, these structures are certain countable pointed directed graphs endowed with a partial vertex labelling with
labels from~$A$. Thus, the number of types is at most $2^\omega\lh A^\omega\le\kappa^\omega$, and consequently, $S_k$ is
$\kappa$-stable whenever $\kappa=\kappa^\omega$.
\end{Pf}

We remark that theories of locally free algebras (including acyclic pairing functions) are also stable; further
model-theoretic properties of acyclic pairing functions were investigated by Bouscaren and Poizat~\cite{bou-poi:pair}.

\section{Computational complexity}\label{sec:comp-compl}

The proof of Theorem~\ref{thm:dec-compl} does not give any bound on the computational complexity of~$S_k$, but as we will
see in this section, we can actually find reasonably tight upper and lower bounds on the complexity. Recall that there
is a general lower bound due to Ferrante and Rackoff~\cite{ferr-rack}:
\begin{Thm}\label{thm:ferr-rack}
Let $T$ be a consistent theory with a pairing function. Then every language $L\in\mathrm{DTIME}(2^0_{O(n)})$ has a
linearly-bounded polynomial-time reduction to~$T$. Consequently, there exists~$\gamma>0$ such that every
decision procedure for~$T$ takes time at least $2^0_{\gamma n}$ for infinitely many input lengths~$n$.
\noproof\end{Thm}

A few remarks are in order. First, the result is stated in~\cite{ferr-rack} for theories of a pairing \emph{function},
but it is straightforward to adapt the argument to theories with a non-functional pairing predicate. The
constant~$\gamma$ only depends on the defining formula for the pairing predicate, otherwise it is independent of~$T$.
Second, the result is quite robust across models of computation and complexity measures: it applies equally well to
time or space, on deterministic, nondeterministic, or alternating Turing machines, etc. The reason is that all these
measures are equivalent up to an exponential or two, and this difference is drowned by the overall complexity: say,
$\mathrm{ASPACE}(2^0_{\gamma n})\sset\mathrm{DTIME}(2^0_{\gamma n+O(1)})$.

In Theorem~\ref{thm:ferr-rack}, the length of input is officially measured as the number of letters when the formula is
written as a word over a finite alphabet (thus a variable $x_i$ takes length $O(\log i)$), but a fortiori the bound
also holds when we measure the input by the number of symbols (quantifiers, connectives, variables, predicate and
function symbols); as we will see, the bound is tight in both regimes (the explanation is that the formulas used in the
lower bound reuse just $O(1)$ distinct variables over and over). We will state upper bounds in terms of the number of
symbols, which is more intuitive, and makes the upper bounds stronger.
\begin{Cor}\label{cor:supext-lb}
There exists $\gamma>0$ such that every decision procedure for any consistent extension of $\vs_2$ has complexity at
least $2^0_{\gamma n}$ for infinitely many input lengths~$n$. In particular, this applies to the theories $S_k$ for
$k\ge2$.
\noproof\end{Cor}

We aim to show that the bound on the complexity of~$S_k$ from Corollary~\ref{cor:supext-lb} is optimal up to the value
of~$\gamma$. The basic idea is that using Theorem~\ref{thm:sim-eleq}, we can represent tuples from an unspecified model
of~$S_k$ by finite objects of bounded size (namely, isomorphism types of $\stcl_m(\ob a)$ for sufficiently large~$m$)
that carry enough information to determine the truth of $\fii(\ob a)$ for a given formula~$\fii$. First, we need an
internal description of structures of the form $\stcl^\sA_m(\ob a)$ so that we can efficiently recognize them.
\begin{Def}\label{def:tcl-str}
Consider a (possibly empty) structure $\sT=\p{T,{\in}^\sT,\ob a}$, where $\Lh(\ob a)=l$. We regard $\sT$ as a directed
graph such that there is an edge $x\to y$ iff $y\in^\sT x$. We say that $\sT$ is a \emph{$\tcl^k_m(l)$-structure} if it
satisfies the following conditions:
\begin{itemize}
\item $\sT$ is a directed acyclic graph with all nodes of out-degree $\le k$.
\item Every node $t\in T$ is reachable from some~$a_i$, $i<l$, in at most $m$ steps.
\item Let $U\sset T$ denote the set of nodes $u\in T$ such that $u$ is reachable from some~$a_i$ in $<m$~steps, or $u$
has out-degree~$k$. Then $\sT$ is extensional w.r.t.~$U$: i.e., for every distinct $u,u'\in U$, there is $t\in T$ such
that $t\in^\sT u$ and $t\notin^\sT u'$, or vice versa.
\end{itemize}
\end{Def}
\begin{Lem}\label{lem:tcl-klm-str}
For any structure $\sT=\p{T,\in^{\sT},\ob a}$ with $\Lh(\ob a)=l$, the following are equivalent:
\begin{enumerate}
\item\label{item:1} $\sT$ is a $\tcl^k_m(l)$-structure.
\item\label{item:2} $\sT$ embeds in a model $\sA\model S_k$ in such a way that $\sT=\stcl^\sA_m(\ob a)$.
\end{enumerate}
\end{Lem}
\begin{Pf}
\ref{item:2}\txto\ref{item:1} is clear, using the observation that $u^\sA\sset\tcl^\sA_m(\ob a)$ for every $u$ from the
set
\[U=\tcl^\sA_{m-1}(\ob a)\cup\bigl\{u\in\tcl^\sA_m(\ob a):\lh{u^\sA\cap\tcl^\sA_m(\ob a)}=k\bigr\}.\]

\ref{item:1}\txto\ref{item:2}: Let $U$ be as in Definition~\ref{def:tcl-str}. We first extend $\sT$ to a
model~$\sA_0=\p{A_0,{\in}^{\sA_0}}$ by adding an infinite descending chain below each $u\in T\bez U$; formally,
$A_0=T\cupd\bigl((T\bez U)\times\omega\bigr)$, with
\[{\in}^{\sA_0}={\in}^\sT\cup\bigl\{\bigl<\p{u,0},u\bigr>,\bigl<\p{u,n+1},\p{u,n}\bigr>:u\in T\bez U,n\in\omega\bigr\}.\]
Notice that no element of~$A_0\bez T$ is reachable in $\le m$ steps from~$\ob a$.
Since each $u\in T\bez U$ has strictly less than~$k$ elements in~$\sT$,
all $u\in A_0$ have at most $k$ elements in~$\sA_0$. Moreover, the structure is still acyclic, and the added chains
ensure that it is extensional; i.e., $\sA_0$ satisfies the axioms $(\ax E)$, $(\ax B_k)$, and $(\ax C_n)$ for all
$n\ge1$.

In order to satisfy axioms $(\ax V_0)$ and~$(\ax V_k)$ as well, we inductively add to~$\sA_0$ all the missing subsets of
size at most~$k$: i.e., we define $\sA_i=\p{A_i,{\in}^{\sA_i}}$ by induction on~$i\in\omega$ as
\begin{align*}
A_{i+1}&=A_i\cupd\{x\sset A_i:\lh x\le k,\forall u\in A_i\:u^{\sA_i}\ne x\},\\
{\in}^{\sA_{i+1}}&={\in}^{\sA_i}\cup\{\p{u,x}:x\in A_{i+1}\bez A_i,u\in x\},
\end{align*}
and we let $\sA=\p{A,{\in}^{\sA}}$ be the union of the chain:
\[A=\bigcup_{i\in\omega}A_i,\qquad{\in}^{\sA}=\bigcup_{i\in\omega}{\in}^{\sA_i}.\]
By construction, $\sA\model S_k$ and $\tcl^\sA_m(\ob a)=\sT$.
\end{Pf}

Incidentally, the construction from Lemma~\ref{lem:tcl-klm-str} can be used to prove Proposition~\ref{prop:fin-ax}:
\begin{Cor}\label{cor:fin-ax}
$S_k$ is not finitely axiomatizable for any $k>0$.
\end{Cor}
\begin{Pf}
Since any finite set of consequences of~$S_k$ is provable from a finite subset of the axiomatization of~$S_k$ in
Definition~\ref{def:sk}, it suffices to show that for every $n\ge1$, there is a model $\sA$ satisfying $(\ax V_0)$,
$(\ax V_k)$, $(\ax E)$, $(\ax B_k)$, $(\ax C_i)$ for $1\le i<n$, and $\neg(\ax C_n)$. Let $\sA_0$ be an $n$-cycle, and
build $\sA$ from~$\sA_0$ as in the proof of Lemma~\ref{lem:tcl-klm-str}.
\end{Pf}
\begin{Def}\label{def:sat-sk}
Let $k,l,m\in\omega$. If $\sT=\p{T,{\in}^\sT,\ob a}$ is a $\tcl^k_m(l)$-structure, and $\fii(\ob x)$ a formula such
that $l=\Lh(\ob x)$ and $m\ge t_k(\rk(\fii))$, we write $\sT\model_{S_k}\fii(\ob a)$ if $\sA\model\fii(\ob a)$, where
$\sA\model S_k$ is such that $\sT=\tcl^\sA_m(\ob a)$. (Such an $\sA$ exists by Lemma~\ref{lem:tcl-klm-str}, and the
definition is independent of the choice of~$\sA$ by Theorem~\ref{thm:sim-eleq}.)

If $\sT=\p{T,{\in}^\sT,\ob a}$ is a $\tcl^k_m(l)$-structure, $m'\le m$, $l'\le l$, and $\ob a'$ is a subsequence of~$\ob a$
of length~$l'$, let $\stcl^\sT_{m'}(\ob a')$ denote the $\tcl^k_{m'}(l')$-structure $\p{T',{\in}^\sT\cap(T'\times
T'),\ob a'}$, where $T'$ is the set of nodes of~$T$ reachable in $\le m'$ steps from~$\ob a'$. (This coincides with
$\stcl^\sA_{m'}(\ob a')$ for any $\sA\model S_k$ such that $\sT=\stcl^\sA_m(\ob a)$.)
\end{Def}

For testing the truth of quantified formulas in~$\sT$, we will need to be able to efficiently recognize when a
$\tcl^k_m(l)$-structure and a $\tcl^k_{m'}(l')$-structure are compatible in that they can be jointly embedded in a
model of~$S_k$. This is accomplished in the next lemma; note that $\ob a$ and~$\ob a'$ need not be disjoint (in fact,
the intended use case is that $\ob a'$ extends $\ob a$).
\begin{Lem}\label{lem:tcl-klm-comp}
Let $\sT=\p{T,{\in}^\sT,\ob a,\ob a'}$, $l=\Lh(\ob a)$, $l'=\Lh(\ob a')$, and $k,m,m'\ge0$. Then the following are
equivalent.
\begin{enumerate}
\item\label{item:7}
$\sT$ embeds in a model $\sA\model S_k$ in such a way that $T=\tcl^\sA_m(\ob a)\cup\tcl^\sA_{m'}(\ob a')$.
\item\label{item:8}
The following conditions hold:
\begin{itemize}
\item $\sT$ is a directed acyclic graph with all nodes of out-degree $\le k$.
\item Every node $t\in T$ is reachable from some~$a_i$, $i<l$, in at most $m$ steps, or from some $a'_i$, $i<l'$, in at
most $m'$ steps.
\item $\sT$ is extensional w.r.t.~$U$, where $U$ denotes the set of nodes $u\in T$ such that $u$ is reachable from
some~$a_i$ in $<m$~steps, or from some~$a'_i$ in $<m'$ steps, or $u$ has out-degree~$k$.
\end{itemize}
\end{enumerate}
\end{Lem}
\begin{Pf}
Just like the proof of Lemma~\ref{lem:tcl-klm-str}.
\end{Pf}

If $\sT$ satisfies the conditions of Lemma~\ref{lem:tcl-klm-comp}, the $\tcl^k_m(l)$-structure $\stcl^\sT_m(\ob a)$ and
the $\tcl^k_{m'}(l')$-structure $\stcl^\sT_{m'}(\ob a')$ are called \emph{compatible}. Note that $\sT$ is uniquely
determined by $\stcl^\sT_m(\ob a)$ and $\stcl^\sT_{m'}(\ob a')$, being their union. We stress that compatibility is
not defined ``up to isomorphism''; the two structures have to be presented in such a way that elements of their
intersection inside~$\sT$ are represented literally the same in both.

We consider the recursive algorithm $\sksat k$ in Fig.~\ref{fig:sk-sat}.
(We are primarily interested in the case where $k$ is a constant, but the algorithm actually works uniformly even if
$k$ is given as part of the input.)
\begin{figure}
\begin{algo}[35em]
\>\key{function} $\sksat k(\sT,\fii)\in\{0,1\}$
\>\key{input:} $\tcl^k_m(l)$-structure $\sT=\p{T,{\in}^\sT,\ob a}$, formula $\fii(\ob x)$,
\>	where $l=\Lh(\ob x)$, $m\ge t_k(\rk(\fii))$
\key{if} $\fii$ is atomic \key{then} \key{return} $\sT\model\fii(\ob a)$
\key{if} $\fii=\neg\fii_0$ \key{then} \key{return} $\neg\sksat k(\sT,\fii_0)$
\key{if} $\fii=\fii_0\lor\fii_1$ \key{then} \key{return} $\sksat k(\sT,\fii_0)\lor\sksat k(\sT,\fii_1)$
\key{if} $\fii=\fii_0\land\fii_1$ \key{then} \key{return} $\sksat k(\sT,\fii_0)\land\sksat k(\sT,\fii_1)$
\key{if} $\fii=\exists y\,\fii_0(\ob x,y)$ \key{then}:
	\key{for each} $\tcl^k_{t_k(\rk(\fii_0))}(l+1)$-structure $\sT'=\p{T',{\in}^{\sT'},\ob a,c}$ \key{do}:
		\key{if} $\sT'$ is compatible with $\sT$ and $\sksat k(\sT',\fii_0)=1$ \key{then} \key{return} 1
	\key{return} 0
\key{if} $\fii=\forall y\,\fii_0(\ob x,y)$ \key{then}:
	\key{for each} $\tcl^k_{t_k(\rk(\fii_0))}(l+1)$-structure $\sT'=\p{T',{\in}^{\sT'},\ob a,c}$ \key{do}:
		\key{if} $\sT'$ is compatible with $\sT$ and $\sksat k(\sT',\fii_0)=0$ \key{then} \key{return} 0
	\key{return} 1
\end{algo}
\caption{An algorithm for $\sT\model_{S_k}\fii(\ob a)$.}
\label{fig:sk-sat}
\end{figure}
\begin{Lem}\label{lem:sk-sat-sound}
Given a $\tcl^k_m(l)$-structure $\sT=\p{T,{\in}^\sT,\ob a}$ and a formula $\fii(\ob x)$ such that $l=\Lh(\ob x)$ and
$m\ge t_k(\rk(\fii))$, $\sksat k(\sT,\fii)=1$ if and only if $\sT\model_{S_k}\fii(\ob a)$.
\end{Lem}
\begin{Pf}
By induction on the complexity of~$\fii$. The only nontrivial cases are for the quantifiers. We will give the proof for
$\fii(\ob x)=\exists y\,\fii_0(\ob x,y)$; the argument for $\forall y\,\fii_0(\ob x,y)$ is dual.

On the one hand, assume that $\sT\model_{S_k}\fii(\ob a)$; i.e., $\sA\model\fii(\ob a)$, where we fix $\sA\model S_k$
such that $\sT=\stcl^\sA_m(\ob a)$. Let $c\in A$ be such that $\sA\model\fii_0(\ob a,c)$, and put
$m'=t_k(\rk(\fii_0))$. Then $\sT'=\stcl^\sA_{m'}(\ob a,c)$ is a $\tcl^k_{m'}(l+1)$-structure compatible with~$\sT$, and
$\sT'\model_{S_k}\fii_0(\ob a,c)$, hence $\sksat k(\sT',\fii_0)=1$ by the induction hypothesis. Thus,
$\sksat k(\sT,\fii)$ returns~$1$ on line~7.

On the other hand, assume that $\sksat k(\sT,\fii)=1$, thus there is a $\tcl^k_{m'}(l+1)$-structure
$\sT'=\p{T',{\in}^{\sT'},\ob a,c}$ compatible with~$\sT$ such that $\sksat k(\sT',\fii_0)=1$. By compatibility, there is
a model $\sA\model S_k$ such that $\sT=\stcl^\sA_m(\ob a)$ and $\sT'=\stcl^\sA_{m'}(\ob a,c)$. By the induction
hypothesis, $\sT'\model_{S_k}\fii_0(\ob a,c)$, which means $\sA\model\fii_0(\ob a,c)$ and $\sA\model\fii(\ob a)$. Thus,
$\sT\model_{S_k}\fii(\ob a)$.
\end{Pf}
\begin{Thm}\label{thm:supexp-ub}
Let $k\ge2$. Given a sentence $\fii$ with $n$~symbols, we can decide whether $S_k\vdash\fii$ in time
$2^{c_k}_{(n+1)/4}$ for sufficiently large~$n$, where $c_k$ is the constant from Proposition~\ref{prop:tk-growth}.
\end{Thm}
\begin{Pf}
We have $S_k\vdash\fii$ iff $\nul\model_{S_k}\fii$ iff $\sksat k(\nul,\fii)=1$ by Lemma~\ref{lem:sk-sat-sound}, where
$\nul$ is considered as a $\tcl^k_m(0)$-structure with $m=t_k(\rk(\fii))$. Rather than measuring time directly, it is
easier to estimate the space requirements of $\sksat k(\nul,\fii)$. We claim that space $O(m\log m+n)$ is sufficient.

It is easy to see that we can test whether a given $\sT$ is a $\tcl^k_m(l)$-structure in space linear in the size of
$\sT$; likewise for testing compatibility, or the truth of atomic formulas. Thus, the dominant cost is that for each
recursive call, we need to store $O(1)$ bits describing where the call was made, and for the quantifier cases, the
structure $\sT'$. The former add up to space $O(n)$, as the recursion depth is at most~$n$. The latter are dominated by
the size of $\sT'$ in the top-most quantifier calls, where it has $s\le k^{\le t_k(\rk(\fii_0))}\le m/k$ elements (in
subsequent calls, the structures become exponentially smaller, hence their space requirements are negligible in
comparison). Since $\sT'$ is a directed graph with out-degree at most~$k$, it can be described by a list of edges using
$O(ks\log s)=O(m\log m)$ bits; this gives total space $O(m\log m+n)$. As long as $m$ dominates~$n$ (which will be the
case for our bounds on~$m$ below), this means the algorithm works in space $O(m\log m)$, and therefore in time
$m^{O(m)}$.

In order to bound $m$ in terms of~$n$, we first bound $r=\rk(\fii)$. Obviously, $r\le n$, but we may do a bit better as
follows. By preprocessing~$\fii$ if necessary, we may assume that there are no dummy quantifiers in~$\fii$. Then each
quantified variable occurs also in an atomic formula; since only two variables occur in a single atomic formula, it
follows that the formula has $\ge r/2$ atomic subformulas (of $3$ symbols each), and consequently $\ge r/2-1$ binary
connectives. Since every quantifier takes two symbols by itself, we see that $n\ge4r-1$, i.e., $r\le(n+1)/4$.

By Proposition~\ref{prop:tk-growth}, $m\le2^{c_k}_{r-1}$, where $r=\rk(\fii)$. Thus, $m\le2^{c_k}_{(n-3)/4}$. Since this grows
much faster than~$n$, we obtain that the algorithm works in space $O(2^{c_k}_{(n-3)/4}\log2^{c_k}_{(n-3)/4})$. In fact,
it is easy to check that there is enough leeway in the bound from Proposition~\ref{prop:tk-growth} so that for any constant~$C$,
$Ct_k(r)\log t_k(r)\le2^{c_k}_{r-1}$ for large enough~$r$. Thus, for large enough~$n$, the algorithm works in space
$2^{c_k}_{(n-3)/4}$, and in time $2^{c_k}_{(n+1)/4}$.
\end{Pf}

The main virtue of Theorem~\ref{thm:supexp-ub} is that it provides an upper bound on the complexity of~$S_k$ that matches the
lower bound from Theorem~\ref{thm:ferr-rack} up to the value of~$\gamma$, and to that end it is stated so that the bound
only depends on~$n$ (and $k$, which is considered to be constant), not other parameters. On the flip side, this
simplicity means that it vastly overestimates the needed complexity for many classes of formulas.

It is clear from the proof that the height of the tower of exponentials in the bound is actually controlled by the
quantifier rank rather than the length of the sentence. Even better, we will show below that it only depends on the
number of quantifier \emph{alternations}.

For simplicity, we will formulate the result for sentences in prenex normal form. Recall that a formula is $\exists_n$
if it is in prenex normal form, and the quantifier prefix consists of $n$ alternating (possibly empty) blocks of
quantifiers, where the first block is existential. The definition of $\forall_n$ formulas is dual. Let us first
generalize Lemma~\ref{lem:simn-ext} and Theorem~\ref{thm:sim-eleq} to handle blocks of quantifiers.
\begin{Lem}\label{lem:min-ext-blk}
Let $\sA,\sB\model S_k$, $\ob a\in A$, $\ob b\in B$, $l=\Lh(\ob a)=\Lh(\ob b)$, and $n,q>0$. If
$\sA,\ob a\sim_{q\,k^{\le n}+n}\sB,\ob b$, then for every $q$-tuple $\ob c\in A$, there exists a $q$-tuple $\ob d\in
B$ such that $\sA,\ob a,\ob c\sim_{n-1}\sB,\ob b,\ob d$.
\end{Lem}
\begin{Pf}
The proof of Lemma~\ref{lem:simn-ext} works literally the same with $\ob c$ in place of~$c$, and $q\,k^{\le n}$ in place of
$k^{\le n}$. In particular, the quantity $k^{\le n}$ only enters the proof through the bound $\lh{\tcl^\sA_n(c)}\le
k^{\le n}$, which is now replaced with $\lh{\tcl^\sA_n(\ob c)}\le q\,k^{\le n}$.
\end{Pf}
\begin{Thm}\label{thm:sim-eleq-blk}
Let $\sA,\sB\model S_k$, $\ob a\in A$, $\ob b\in B$, and $\Lh(\ob a)=\Lh(\ob b)$. For any $n,q\in\omega$, define
$t_k(n,q)$ by $t_k(0,q)=0$, $t_k(n+1,q)=q\,k^{\le t_k(n,q)+1}+t_k(n,q)+1$. Let $\fii(\ob x)$ be an $\exists_n$ formula
with each quantifier block of length at most~$q$. Then
\[\sA,\ob a\sim_{t_k(n,q)}\sB,\ob b\implies\bigl(\sA\model\fii(\ob a)\iff\sB\model\fii(\ob b)\bigr).\]
\end{Thm}
\begin{Pf}
By induction on~$n$, using Lemma~\ref{lem:min-ext-blk}.
\end{Pf}
\begin{Lem}\label{lem:tk-growth-blk}
For any $k\ge2$ and $n,q\ge1$, we have
\[t_k(n,q)\le2^{(q(k+1)+2)\log k+\log\log k+\log q+2}_{n-1}\le2^{4qk\log k}_{n-1}.\]
\end{Lem}
\begin{Pf}
Similar to the proof of Proposition~\ref{prop:tk-growth}, with $(x+1)\log k+\log\log k+\log q+2$ in place of $h(x)$, using the
inequality
\[(q\,k^{\le x+1}+x+2)\log k+\log\log k+\log q+2\le4q\,k^{x+1}\log k,\]
which can be proved in the same way as~\eqref{eq:7}.
\end{Pf}
\begin{Thm}\label{thm:qalt-ub}
Given a sentence $\fii$ in prenex normal form and $k\ge2$, we can decide whether $S_k\vdash\fii$ in
$\cxt{NTIME}\bigl(t_k(r,q)^{O(t_k(r,q))}n^{O(1)}\bigr)$, where $n$ is the length of~$\fii$, $r$ is such that
$\fii$ is $\exists_{r+1}$, and $q$~is the maximal length of a quantifier block in~$\fii$. This is
$\cxt{NTIME}(n^{O(1)})$ for $r=0$, $\cxt{NTIME}\bigl((kq)^{O(kq)}n^{O(1)}\bigr)$ for $r=1$, and
$\cxt{NTIME}\bigl(2^{O(qk\log k)}_rn^{O(1)}\bigr)$ for $r\ge2$.
\end{Thm}
\begin{figure}
\begin{algo}[35em]
\>\key{function} $\bsksat k(\sT,\fii)\in\{0,1\}$
\>\key{input:} $\tcl^k_m(l)$-structure $\sT=\p{T,{\in}^\sT,\ob a}$, $\exists_r$ or $\forall_r$ formula $\fii(\ob x)$,
\>	where $l=\Lh(\ob x)$, $m\ge t_k(r,q)$, $q=$ maximal quantifier block size in $\fii$
\key{if} $r=0$ \key{then} \key{return} $\sT\model\fii(\ob a)$
\key{if} $\fii=\exists\ob y\,\fii_0(\ob x,\ob y)$, $\fii_0\in\forall_{r-1}$, $l_0=\Lh(\ob y)$ \key{then}:
	\key{for each} $\tcl^k_{t_k(r-1,q)}(l+l_0)$-structure $\sT'=\p{T',{\in}^{\sT'},\ob a,\ob c}$ \key{do}:
		\key{if} $\sT'$ is compatible with $\sT$ and $\bsksat k(\sT',\fii_0)=1$ \key{then} \key{return} 1
	\key{return} 0
\key{if} $\fii=\forall\ob y\,\fii_0(\ob x,\ob y)$, $\fii_0\in\exists_{r-1}$, $l_0=\Lh(\ob y)$ \key{then}:
	\key{for each} $\tcl^k_{t_k(r-1,q)}(l+l_0)$-structure $\sT'=\p{T',{\in}^{\sT'},\ob a,\ob c}$ \key{do}:
		\key{if} $\sT'$ is compatible with $\sT$ and $\bsksat k(\sT',\fii_0)=0$ \key{then} \key{return} 0
	\key{return} 1
\end{algo}
\caption{A block-wise algorithm for $\sT\model_{S_k}\fii(\ob a)$.}
\label{fig:sk-sat-blk}
\end{figure}
\begin{Pf}
Write $\fii=\exists\ob x\,\psi(\ob x)$, where $\psi$ is $\forall_r$. Put $m=t_k(r,q)$ and $l=\Lh(\ob x)\le q$. In order
to test $S_k\vdash\fii$, we nondeterministically guess a $\tcl^k_m(l)$-structure $\sT=\p{T,{\in}^\sT,\ob a}$, and
verify $\sT\model_{S_k}\psi(\ob a)$ using the algorithm $\bsksat k(\sT,\psi)$ from Fig.~\ref{fig:sk-sat-blk}. Note that
$s=\lh T\le l\,k^{\le m}\le\frac1kt_k(r+1,q)$, hence the bit-size of~$\sT$ is
$O(ks\log s)=O\bigl(t_k(r+1,q)\log t_k(r+1,q)\bigr)$, and we can check that $\sT$ is a $\tcl^k_m(l)$-structure in time
polynomial in $t_k(r+1,q)=q\,k^{O(t_k(r,q))}$. (For $r=0$, we have $s\le l\le n$, thus $\sT$ can be represented with
$O(n^2)$ bits using an adjacency matrix, and then we can check that $\sT$ is a $\tcl^k_0(l)$-structure in time
$n^{O(1)}$ independent of~$k$: if $k>n$, we only need to check that $T=\{\ob a\}$ and $\in^\sT$ is acyclic.)

We claim that $\bsksat k(\sT,\psi)$, and thus the whole test, works in time polynomial in $n$ and $t_k(r,q)^{t_k(r,q)}$.
In the top-level iteration, the structures $\sT'$ have sizes up to
$(l+l_0)k^{\le t_k(r-1,q)}\le2q\,k^{\le t_k(r-1,q)}\le\frac2kt_k(r,q)$, and can be described using
$O\bigl(t_k(r,q)\log t_k(r,q)\bigr)$ bits. Thus, the loop on lines 7--8 goes through
$\exp\bigl(O\bigl(t_k(r,q)\log t_k(r,q)\bigr)\bigr)$ structures $\sT'$; for each of them, it checks in time
$t_k(r+1,q)^{O(1)}$ whether it is compatible with~$\sT$, and if so, makes a recursive call. In turn, each of these
recursive calls will involve a loop over $\exp\bigl(O\bigl(t_k(r-1,q)\log t_k(r-1,q)\bigr)\bigr)$ structures,
where for each of them, we do a compatibility check in time $t_k(r,q)^{O(1)}$, and a recursive call. This goes on until
we get down to the quantifier-free matrix at recursion depth~$r$; this takes time $n^{O(1)}$ to check on line~1. Thus, the total
number of recursive calls is
\begin{align*}
\prod_{i<r}2^{O(t_k(r-i,q)\log t_k(r-i,q))}&=2^{O\bigl(\sum_{i<r}t_k(r-i,q)\log t_k(r-i,q)\bigr)}\\
&=2^{O(t_k(r,q)\log t_k(r,q))}=t_k(r,q)^{O(t_k(r,q))},
\end{align*}
and each takes time polynomial in $t_k(r+1,q)=q\,k^{O(t_k(r,q))}$ and~$n$. This gives total time
$t_k(r,q)^{O(t_k(r,q))}n^{O(1)}$, as claimed. (For $r=0$, this means $n^{O(1)}$; there are no
recursive calls.)

For $r=1$, we have $t_k(1,q)=O(kq)$, hence the time bound is $(kq)^{O(kq)}n^{O(1)}$. For $r\ge2$, $t_k(r,q)=2^{O(qk\log
k)}_{r-1}$ by Lemma~\ref{lem:tk-growth-blk}; it is easy to show that $(2^x_d)^c\le2^{cx}_d$ for any $d,x\ge1$ by induction
on $d$, hence
\[t_k(r,q)^{O(t_k(r,q))}=2^{t_k(r,q)^{O(1)}}=2^{\bigl(2^{O(qk\log k)}_{r-1}\bigr)^{O(1)}}
  =2^{2^{O(qk\log k)}_{r-1}}=2^{O(qk\log k)}_r,\]
which gives the time bound $2^{O(qk\log k)}_rn^{O(1)}$.
\end{Pf}

For completeness, let us also indicate the complexity of $S_k$ for $k=0,1$, which is essentially known from the
literature.
\begin{Thm}\label{thm:cxt-s0-s1}
$S_1$ is $\psp$-complete, and for any fixed $r\ge1$, the $\exists_r$~fragment of~$S_1$ is $\Sigma^\ptime_r$-complete,
and the $\forall_r$~fragment is $\Pi^\ptime_r$-complete.

$S_0$ is decidable in~$\ptime$; more precisely, it is complete for $\cxt{ALOGTIME}=U_E\text{-uniform }\nci$ under
$\LT$ reductions.
\end{Thm}
\begin{Pf}
$S_0$ is the theory of a one-element structure, hence it is equivalent to propositional logic (we can decide a given
sentence by removing all quantifiers, replacing atomic formulas $x_i\in x_j$ with the truth-constant~$0$ and
$x_i=x_j$ with~$1$, and evaluating the resulting Boolean sentence). This is $\nci$-complete by results of
Buss~\cite{buss:bfvp}.

For $S_1$, it is well known and easy to see that the truth of quantified Boolean sentences is reducible to any
consistent first-order theory~$T$ that proves the existence of two distinct elements, making $T$ $\psp$-hard. Moreover,
the reduction takes $\exists_r$~QBF to $\exists_r$~sentences, hence $T$-provability of $\exists_r$~sentences is
$\Sigma^\ptime_r$-hard, and dually for $\forall_r$.

On the other hand, $\sH_1$ is definitionally equivalent to $\p{\N,0,S}$, whose theory is known to be decidable in
$\psp$: e.g., this is proved in~\cite{ferr-rack} for the more general structure $\p{\N,{<}}$.

Using the machinery we have already developed, this can be shown as follows. First, we have
$t_1(r,q)=(2+\frac1q)(q+1)^r$, hence $\bsksat1(\nul,\fii)$ decides $S_1\vdash\fii$ in exponential space. We can make it
more space-efficient by employing a more compact representation for $\tcl^1_m(l)$-structures
$\sT=\p{T,{\in}^\sT,\ob a}$. Any such structure is a disjoint union of $\in^\sT$-chains, where each chain has some
$a_i$ on top, the distance between neighbouring $a_i$ and~$a_j$ on the same chain is $\le m+1$, and each chain ends
$\le m$ steps below the lowest $a_i$ on the chain. We can represent this by noting for each $a_i$ the nearest $a_j$
below~$a_i$ on the same chain, if any, and the distance (in binary!) from $a_i$ to~$a_j$, or to the end of the chain.
This takes $O(l\log m)=O(lr\log(q+1))=O(n^2\log n)$ bits if $m=t_1(r,q)$ and $l,r,q\le n$, and we can test
compatibility of $\tcl^1_m(l)$-structures in this representation and satisfaction of quantifier-free formulas in
polynomial time.

Thus, $\bsksat1$ modified to use this representation runs in polynomial space, placing $S_1$ in $\psp$. Moreover, we
may view the modified $\bsksat1$ as an alternating polynomial-time algorithm, where the loop on lines 3--4 is replaced
with a nondeterministic (existential) guess of $\sT'$, and the loop on lines 7--8 with a co-nondeterministic
(universal) guess. Then $\bsksat1(\nul,\fii)$ for an $\exists_r$~sentence $\fii$ makes $r-1$ alternations starting from
an existential state, i.e., it works in $\Sigma_r\text-\cxt{TIME}(n^{O(1)})=\Sigma^\ptime_r$, and dually for
$\forall_r$~sentences.
\end{Pf}

\section{Conclusion}\label{sec:conclusion}

As we have seen, the complete theory of the structure $\sH_k$ can be described by a simple list of axioms, it is
decidable, and generally tame (it has quantifier elimination down to formulas of quite a low complexity, it is
stable, and its computational complexity---albeit somewhat daunting---is the lowest possible for theories with pairing).
Thus, it is in many respects as nice as other known examples of decidable theories with a pairing function.

However, it has a different flavour from the previous examples, which are generally of algebraic or arithmetic nature,
whereas here we have a theory of sets. In particular, the theories $\Th(\sH_k)$ provide natural decidable extensions of
finite fragments of the Vaught set theory~$\vs$, which was our original motivation.

\section*{Acknowledgements}
I want to thank Albert Visser for a helpful discussion on theories of pairing, and an anonymous reviewer for useful
comments. A preliminary sketch of the basic results of this paper was first reported in~\cite{ej:hk-mo}.

\bibliographystyle{mybib}
\bibliography{hered-bounded}

\ifx\url\undefined {\catcode`\/=13
  \gdef/{\string/\futurelet\nexttoken\finishslash}
  \gdef\finishslash{\ifx\nexttoken/\else\penalty\relpenalty\fi}}
  \def\url{\begingroup\catcode`\~=12 \catcode`\_=12 \catcode`\/=13 \finishurl}
  \def\finishurl#1{\texttt{#1}\endgroup} \fi
\providecommand{\bysame}{\leavevmode\hbox to5em{\hrulefill}\thinspace}
\providecommand\bibliographyhook{}
\begin{thebibliography}{10}
\bibliographyhook

\bibitem{ack:hered-fin}
Wilhelm Ackermann, \emph{Die {W}iderspruchsfreiheit der allgemeinen
  {M}engenlehre}, Mathematische Annalen 114 (1937), pp.~305--315 (in German).

\bibitem{bou-poi:pair}
Elisabeth Bouscaren and Bruno Poizat, \emph{Des belles paires aux beaux uples},
  Journal of Symbolic Logic 53 (1988), no.~2, pp.~434--442 (in French).

\bibitem{buss:bfvp}
Samuel~R. Buss, \emph{The {B}oolean formula value problem is in
  {$\mathit{ALOGTIME}$}}, in: Proceedings of the 19th Annual {ACM} {S}ymposium
  on {T}heory of {C}omputing, ACM Press, 1987, pp.~123--131.

\bibitem{ceg-gri-rich:can-pair}
Patrick C{\'e}gielski, Serge Grigorieff, and Denis Richard, \emph{La
  th{\'e}orie {\'e}l{\'e}mentaire de la fonction de couplage de {C}antor des
  entiers naturels est d{\'e}cidable}, Comptes Rendus de l'Acad{\'e}mie des
  Sciences -- Series I -- Mathematics 331 (2000), no.~2, pp.~107--110 (in
  French).

\bibitem{ceg-rich:lists}
Patrick C{\'e}gielski and Denis Richard, \emph{On arithmetical first-order
  theories allowing encoding and decoding of lists}, Theoretical Computer
  Science 222 (1999), no.~1--2, pp.~55--75.

\bibitem{ceg-rich:can-pair-succ}
\bysame, \emph{Decidability of the theory of the natural integers with the
  {C}antor pairing function and the successor}, Theoretical Computer Science
  257 (2001), no.~1--2, pp.~51--77.

\bibitem{ferr-rack}
Jeanne Ferrante and Charles~W. Rackoff, \emph{The computational complexity of
  logical theories}, Lecture Notes in Mathematics vol. 718, Springer-Verlag,
  1979.

\bibitem{ej:hk-mo}
Emil Je{\v r}{\'a}bek, \emph{Answer to a question by {Z}uhair {A}l-{J}ohar},
  MathOverflow, 2018, \url{https://mathoverflow.net/a/302634}.

\bibitem{malc:loc-free}
Anatoli{\u\i}~I. Mal'cev, \emph{On the elementary theories of locally free
  universal algebras}, Doklady Akademii Nauk SSSR 138 (1961), no.~5,
  pp.~1009--1012 (in Russian), {E}nglish translation in: {S}oviet {M}athematics
  -- Doklady 2 (1961), no.~3, pp.~768--771.

\bibitem{malc:ax-cl-loc-free}
\bysame, \emph{Axiomatizable classes of locally free algebras of several
  types}, Sibirski{\u\i} Matematicheski{\u\i} Zhurnal 3 (1962), no.~5,
  pp.~729--743 (in Russian), {E}nglish translation in: The Metamathematics of
  Algebraic Systems: Collected Papers: 1936--1967, Studies in Logic and the
  Foundations of Mathematics 66, North-Holland, 1971, pp.~262--281.

\bibitem{sem:plus-pow}
Aleksei~L. Semenov, \emph{Logical theories of one-place functions on the set of
  natural numbers}, Izvestiya Akademii Nauk SSSR, Seriya Matematicheskaya 47
  (1983), no.~3, pp.~623--658 (in Russian), {E}nglish translation in:
  Mathematics of the USSR, Izvestiya 22 (1984), no.~3, pp.~587--618.

\bibitem{tmr}
Alfred Tarski, Andrzej Mostowski, and Rafael~M. Robinson, \emph{Undecidable
  theories}, North-Holland, Amsterdam, 1953.

\bibitem{tenn:phd}
Richard~L. Tenney, \emph{Decidable pairing functions}, Ph.D. thesis, Cornell
  University, 1972.

\bibitem{vaught}
Robert~L. Vaught, \emph{Axiomatizability by a schema}, Journal of Symbolic
  Logic 32 (1967), no.~4, pp.~473--479.

\bibitem{vis:pair}
Albert Visser, \emph{Pairs, sets and sequences in first-order theories},
  Archive for Mathematical Logic 47 (2008), no.~4, pp.~299--326.

\end{thebibliography}
\end{document}